# Differentiation and Ordered Optimization in Banach Spaces

Jinlu Li

Department of Mathematics
Shawnee State University
Portsmouth, Ohio 45662 USA
jli@shawnee.edu

**Abstract.** This paper is devoted to extend the applications of differentiation to some optimization problems in ordinary calculus to the applications of Gâteaux and Fréchet differentiation to some ordered optimization-related problems in ordered Banach spaces. In order to achieve this purpose, we will introduce the concepts of generalized critical point, ordered extreme and order monotone property of single-valued mappings in partially ordered Banach spaces, which are extension of critical points in ordinary calculus. Then, we will provide the connections between generalized critical point, ordered extreme and Gâteaux and Fréchet derivatives of single-valued mappings in partially ordered Banach spaces that analogously extends the connection between critical points and extrema of real valued functions in calculus. In contrast with the explicit formulas of derivatives in calculus, we will find the explicit formulas of Gâteaux and Fréchet derivatives of some single-valued mappings on the Banach spaces $(l_p, \|\cdot\|_p)$, for $1 < p < \infty$ and $C[0, 1]$, such as polynomial type operators and trigonometric type operators. By these concepts, we will investigate the connection between generalized critical points and ordered extrema of single-valued mappings in partially ordered Banach spaces that extends the connection between critical points and extrema of real valued functions in calculus. In the last 3 sections, we will prove that in partially ordered Banach spaces, the order monotone of single-valued mappings can be described by its Gâteaux derivatives or Fréchet derivatives in partially ordered Banach spaces.

*Dedicate to Professor George Xianzhi Yuan in the Occasion of His 60th Birthday*

## 1. Introduction and Preliminary

In modern analysis, as a natural extension of the classical differentiability of functions in calculus, the differentiability of mappings between normed vector spaces (in particular in Banach spaces) has been developed rapidly. The most popular and most useful concepts of differentiation in Banach spaces are Gâteaux (directional) differentiability and Fréchet differentiability. For reader's convenience of reference, we review these concepts below.

Throughout this paper, unless otherwise stated, let $(X, \|\cdot\|_X)$ and $(Y, \|\cdot\|_Y)$ be Banach spaces with origins $\theta_X$ and $\theta_Y$, respectively. Let $A$ be a nonempty convex and open subset of $X$ and let $\bar{x} \in A$. Let $T: A \to Y$ be a single-valued mapping.

(Gâteaux (directional) differentiability of $T$). Let $v \in X$ with $v \neq \theta_X$. If there is a vector in $Y$ denoted by $T'(\bar{x}, v)$ such that

$$T'(\bar{x}, v) = \lim_{t\to 0}\frac{T(\bar{x}+tv)-T(\bar{x})}{t}, \tag{1.1}$$

then, $T$ is said to be Gâteaux directionally differentiable at point $\bar{x}$ along direction $v$. The point $T'(\bar{x}, v) \in Y$ is called the Gâteaux directional derivative of $T$ at point $\bar{x}$ along direction $v$. Furthermore, if $T$ is Gâteaux directionally differentiable at point $\bar{x}$ along every direction $v \in X\backslash\{\theta_X\}$, then $T$ is said to be Gâteaux differentiable at $\bar{x}$ and the Gâteaux derivative of $T$ at point $\bar{x}$ is denoted by

$$T'(\bar{x})(v) = \lim_{t \to 0} \frac{T(\bar{x}+tv)-T(\bar{x})}{t}, \text{ for any } v \in X\backslash\{\theta_X\}. \tag{1.2}$$

(Fréchet differentiability of $T$). If there is a continuous and linear mapping $\nabla T(\bar{x})$: $X \rightarrow Y$ such that,

$$\lim_{u \xrightarrow{X} \theta_X} \frac{T(\bar{x}+u)-T(\bar{x})-\nabla T(\bar{x})(u)}{\|u\|_X} = \theta_Y, \tag{1.3}$$

then $T$ is said to be Fréchet differentiable at $\bar{x}$ and $\nabla T(\bar{x})$ is called the Fréchet derivative of $T$ at $\bar{x}$. It is well-known that in Banach spaces (or in normed vector spaces), for single-valued mappings, Fréchet differentiability is stronger than Gâteaux differentiability. That is,

$$T \text{ is Fréchet differentiable at } \bar{x} \quad \Longrightarrow \quad T \text{ is Gâteaux differentiable at } \bar{x}. \tag{1.4}$$

More precisely, the connection between Gâteaux differentiability and Fréchet differentiability in Banach spaces is that if $T$ is Fréchet differentiable at $\bar{x}$, then $T$ is also Gâteaux differentiable at $\bar{x}$ and

$$T'(\bar{x})(v) = \nabla T(\bar{x})(v), \text{ for any } v \in X\backslash\{\theta_X\}. \tag{1.5}$$

Both Gâteaux and Fréchet derivatives have many properties, which include linearity property. Let $T_1$, $T_2$: $A \rightarrow Y$ be single-valued mappings. Let $\bar{x} \in A$ and $v \in X\backslash\{\theta_X\}$ and let $a_1$, $a_2$ be scalars. Suppose that both $T_1$ and $T_2$ are Gâteaux differentiable at $\bar{x}$. Then $a_1T_1 + a_2T_2$ is also Gâteaux differentiable at $\bar{x}$ and

$$(a_1T_1 + a_2T_2)'(\bar{x})(v) = a_1T_1'(\bar{x})(v) + a_2T_2'(\bar{x})(v), \text{ for any } v \in X\backslash\{\theta_X\}. \tag{1.6}$$

If both $T_1$ and $T_2$ are Fréchet differentiable at $\bar{x}$. Then $a_1T_1 + a_2T_2$ is also Fréchet differentiable at $\bar{x}$ and

$$\nabla(a_1T_1 + a_2T_2)(\bar{x})(v) = a_1\nabla T_1(\bar{x})(v) + a_2\nabla T_2(\bar{x})(v), \text{ for any } v \in X\backslash\{\theta_X\}. \tag{1.7}$$

In contrast with ordinary derivatives in calculus, Fréchet derivatives also satisfy the following chain rule in Banach spaces.

Let $(Z, \|\cdot\|_Z)$ be a Banach space. Let $A$ and $E$ be nonempty convex open subsets in $X$ and $Y$, respectively. Let $T$: $A \rightarrow E$ and $S$: $E \rightarrow Z$ be single-valued mappings. Let $\bar{x} \in A$ and $\bar{y} \in E$ with $\bar{y} = T(\bar{x})$. Suppose that $T$ is Fréchet differentiable at $\bar{x}$ and $S$ is Fréchet differentiable at $\bar{y}$. Then the composite single valued mapping $S \circ T$: $A \rightarrow Z$ is Fréchet differentiable at $\bar{x}$ and

$$\nabla(S \circ T)(\bar{x}) = \nabla S(\bar{y}) \circ \nabla T(\bar{x}). \tag{1.8}$$

Both Gâteaux (directional) differentiation and Fréchet differentiation in Banach spaces have been studied by many authors (See [1, 4, 6, 10, 11, 16, 19−24]). Furthermore, these concepts have been generalized to general topological vector spaces with some applications (See [27, 33]). Gâteaux and Fréchet derivatives of single-valued mappings in Banach spaces have been widely applied to many fields in both pure and applied mathematics (See [9, 12, 13]), in particular in optimization theory (See [2, 3, 7, 8, 14, 15, 17, 28−31]), variational analysis (See [5,18, 25, 26, 32, 34]), game theory (See [2, 3, 28−31]), and so forth.

Notice that the monotone properties of real valued functions in calculus are defined with respect to the standard order of real numbers. For the purposes of dealing with ordered optimizations in Banach spaces, as usual, we will consider the underlying Banach spaces $(X, \|\cdot\|_X)$ and $(Y, \|\cdot\|_Y)$ to be equipped with partial orders $\preccurlyeq_C$ and $\preccurlyeq_K$, respectively, where, $\preccurlyeq_C$ and $\preccurlyeq_K$ are induced by nontrivial closed convex and pointed cones $C$ and $K$ in $X$ and $Y$, respectively. With respect to the partial orders equipped on Banach spaces, we can define ordered extrema and order monotone for single-valued mappings between partially ordered Banach spaces.

In sections 4 and 5, we will use Gâteaux and Fréchet derivatives to define generalized critical points of single-valued mappings between Banach spaces, which extend the concept of critical points of real valued functions in calculus. If $T$ is Gâteaux differentiable at $\bar{x}$ and $T'(\bar{x})(v) = \theta_Y$, for every $v \in X\backslash\{\theta_X\}$, then $\bar{x}$ is called a generalized critical point of $T$. If for $x_1, x_2 \in A$, we have that $x_1 \preccurlyeq_C x_2$ implies that the values of $T$ will satisfies $T(x_1) \preccurlyeq_K T(x_2)$, then, $T$ is said to be $\preccurlyeq_C$-$\preccurlyeq_K$ increasing on $A$. For notational simplicity, $T$ is said to be order increasing on $A$. The concept of order decreasing on $A$ can be similarly defined.

In contrast with the monotone properties and extrema in calculus, in section 4, we will show that, if $T$ is Gâteaux differentiable at $\bar{x}$, then

$$\bar{x} \text{ is an } \preccurlyeq_K\text{-extreme point of } T \quad \Longrightarrow \quad \bar{x} \text{ is a generalized critical point of } T. \tag{1.9}$$

Similar to the property of extrema in calculus, by using the Banach space $l_p$, we will provide some counter examples to show that

$$\bar{x} \text{ is a generalized critical point of } T \quad \not\Rightarrow \quad \bar{x} \text{ is an } \preccurlyeq_K\text{-extreme point of } T. \tag{1.10}$$

To extend monotone properties to partially ordered Banach spaces, in section 5, we will prove the connection between the order monotone and Gâteaux and Fréchet derivatives as follows.

$$T \text{ is order increasing on } A \quad \Longrightarrow \quad T'(x)(v) \succcurlyeq_K \theta_Y, \text{ for any } x \in A \text{ and for any } v \in C\backslash\{\theta_X\}. \tag{1.11}$$

We will give some applications of Gâteaux and Fréchet derivatives of single valued mappings to ordered optimization problems between partially ordered Banach spaces.

## 2. Generalized Critical Points and Ordered Extrema in Partially Ordered Banach Space

### 2.1. Generalized Critical Point and Ordered Extrema in Banach Space

Throughout this section, as stated in section 1, we let $(X, \|\cdot\|_X)$ and $(Y, \|\cdot\|_Y)$ be Banach spaces with origins $\theta_X$ and $\theta_Y$. For given $x_0 \in X$ and $r > 0$, let $B_X(x_0, r)$ denote the closed ball in $X$ with radius $r$ and centered at $x_0$. For given $y_0 \in Y$ and $r > 0$, we similarly define the closed ball $B_Y(y_0, r)$ in $Y$ with radius $r$ and centered at $y_0$. Let $A$ be a nonempty convex open subset in $X$. Let $T: A \to Y$ be a single-valued mapping and let $\bar{x} \in A$. We review the following definitions from [27].

**Definition 2.1.** Suppose that $T$ is Gâteaux differentiable at $\bar{x}$. If

$$T'(\bar{x})(v) = \theta_Y, \text{ for every } v \in X\backslash\{\theta_X\},$$

then $\bar{x}$ is called a generalized critical point of $T$.

Let $\theta_{X,Y}$ denote the null continuous and linear operator from $X$ to $Y$ with $\theta_{X,Y}(x) = \theta_Y$, for any $x \in X$.

**Lemma 2.2**. *Suppose that $T$ is Fréchet differentiable at $\bar{x}$. If $\nabla(\bar{x}) = \theta_{X,Y}$, then $\bar{x}$ is a generalized critical point of $T$.*

*Proof.* This lemma is proved by Definition 2.1 and the connection between Gâteaux derivatives and

Fréchet derivatives in Banach spaces that is presented in (1.5). □

We write the following sets of generalized critical points of $T$:

$$dT^0 \coloneqq \{\bar{x} \in A: T'(\bar{x}) \text{ exists and } T'(\bar{x})(v) = \theta_Y, \text{for any } v \in X\backslash\{\theta_X\}\};$$

$$\nabla T^0 \coloneqq \{\bar{x} \in A: \nabla(\bar{x}) \text{ exists and } \nabla(\bar{x}) = \theta_{X,Y}\}.$$

$dT^0$ is the set of generalized critical point of $T$. By Lemma 2.2, we have

$$\nabla T^0 \subseteq dT^0. \tag{2.1}$$

**Definition 2.3.** Let $\bar{x} \in A$ and let $v \in X\backslash\{\theta_X\}$. If

$$T(\bar{x} + tv) \preccurlyeq_K T(\bar{x}), \text{ for all real } t \text{ with } \bar{x} + tv \in A,$$

then, $T$ is said to have an $\preccurlyeq_K$-maximum value in $A$ at $\bar{x}$ along direction $v$. $T(\bar{x})$ is called the $\preccurlyeq_K$-maximum value of $T$ in $A$ along direction $v$. The $\preccurlyeq_K$-minimum value of $T$ at $\bar{x}$ along direction $v \in X\backslash\{\theta_X\}$ can be similarly defined. More strictly, if $T$ satisfies the following order-inequality

$$T(x) \preccurlyeq_K T(\bar{x}), \text{ for all } x \in A, \tag{2.2}$$

then, $T$ is said to have an (absolute) $\preccurlyeq_K$-maximum value in $A$ at $\bar{x}$ and $T(\bar{x})$ is called the (absolute) $\preccurlyeq_K$-maximum value of $T$ in $A$. The (absolute) $\preccurlyeq_K$-minimum value of $T$ in $A$ at $\bar{x}$ can be similarly defined. In either of (absolute) $\preccurlyeq_K$-maximum value or (absolute) $\preccurlyeq_K$-minimum value, they are called $\preccurlyeq_K$-extreme values. In this case, we say that $\bar{x}$ is an $\preccurlyeq_K$-extreme point of the mapping $T$. With respect to the given set $A$ and the mapping $T: A \to Y$, we write the sets of $\preccurlyeq_K$-extreme points of the mapping $T$.

$$E(T, \preccurlyeq_K) = \{\bar{x} \in A: \bar{x} \text{ is an } \preccurlyeq_K\text{-extreme point of } T\};$$
$$E^+(T, \preccurlyeq_K) = \{\bar{x} \in A: \bar{x} \text{ is an } \preccurlyeq_K\text{-maximum point of } T\};$$
$$E^-(T, \preccurlyeq_K) = \{\bar{x} \in A: \bar{x} \text{ is an } \preccurlyeq_K\text{-minimum point of } T\}.$$

If the partial order $\preccurlyeq_K$ is understand in the considered problem and it is not confused, then $E(T, \preccurlyeq_K)$, $E^+(T, \preccurlyeq_K)$ and $E^-(T, \preccurlyeq_K)$ are simply denoted by $E(T)$, $E^-(T)$ and $E^-(T)$, respectively.

In the next theorem, we investigate the connection between generalized critical points and ordered extreme points of single-valued mappings between partially ordered Banach spaces, which naturally generalizes the connection between critical points and extreme points of real valued functions in calculus. In general, this theorem has been proved in [27] for partially ordered Hausdorff topological vector spaces. Here, we will reprove it for partially ordered Banach spaces with proof slightly different from the proof of Theorem 9.4 in [27].

**Theorem 2.4**. *Let $\bar{x} \in A$ and $v \in X\backslash\{\theta_X\}$. If $T$ has $\preccurlyeq_K$-maximum ($\preccurlyeq_K$-minimum) value at $\bar{x}$ in direction $v$ and $T$ is Gâteaux directionally differentiable at $\bar{x}$ in direction $v$, then $T'(\bar{x}, v) = \theta_Y$.*

*Proof.* Suppose that $T$ is Gâteaux directionally differentiable at point $\bar{x}$ along the given (fixed) direction $v$ and $T$ takes an $\preccurlyeq_K$-minimum value at point $\bar{x}$ along direction $v$ (the $\preccurlyeq_K$-maximum case can be similarly proved). At first, we prove $T'(\bar{x}, v) \succcurlyeq_K \theta_Y$, this is equivalent to $T'(\bar{x}, v) \in K$. Assume, by the way of contradiction, that $T'(\bar{x}, v) \notin K$, this implies $T'(\bar{x}, v) \in Y\backslash K$. Since $Y\backslash K$ is an open subset of $Y$, by $T'(\bar{x}, v) \in Y\backslash K$, there is $\varepsilon > 0$ such that, $B_Y^o(T'(\bar{x}, v), \varepsilon) \subseteq Y\backslash K$, where $B_Y^o(T'(\bar{x}, v), \varepsilon)$ is the interior of $B_Y(T'(\bar{x}, v), \varepsilon)$ satisfying

$$B_Y^o(T'(\bar{x}, v), \varepsilon) = \{y \in Y : \|y - T'(\bar{x}, v)\|_Y < \varepsilon\}. \tag{2.3}$$

The condition $B_Y^o(T'(\bar{x}, v), \varepsilon) \subseteq Y \setminus K$ induces

$$B_Y^o(T'(\bar{x}, v), \varepsilon) \cap K = \emptyset. \tag{2.4}$$

(2.3) and (2.4) together imply that, for $y \in Y$, we have

$$y \in K \quad \Longrightarrow \quad \|y - T'(\bar{x}, v)\|_Y \geq \varepsilon. \tag{2.5}$$

Since $T$ is Gâteaux directionally differentiable at point $\bar{x}$ along the given (fixed) direction $v$, by definition and by the condition that $A$ is an open subset in $X$ and $\bar{x} \in A$, for the above given $\varepsilon > 0$, there is $\delta > 0$ such that, for real $t$,

$$0 < |t| < \delta \quad \Longrightarrow \bar{x} + tv \in A \ \text{ and } \ \left\|\frac{T(\bar{x}+tv)-T(\bar{x})}{t} - T'(\bar{x}, v)\right\|_Y < \varepsilon. \tag{2.6}$$

In particular, let $t > 0$ in (2.6), we get

$$0 < t < \delta \quad \Longrightarrow \quad \bar{x} + tv \in A \ \text{ and } \ \left\|\frac{T(\bar{x}+tv)-T(\bar{x})}{t} - T'(\bar{x}, v)\right\|_Y < \varepsilon. \tag{2.7}$$

By the assumption that $T(\bar{x} + tv) - T(\bar{x}) \succcurlyeq_K \theta_Y$, we must have that

$$0 < t < \delta \quad \Longrightarrow \quad \frac{T(\bar{x}+tv)-T(\bar{x})}{t} \in K. \tag{2.8}$$

For $t > 0$ in (2.8) and (2.7), this creates a contradiction between (2.7) and (2.5). Hence $T'(\bar{x}, v) \in K$ and

$$T'(\bar{x}, v) \succcurlyeq_K \theta_Y. \tag{2.9}$$

Next, we prove $T'(\bar{x}, v) \preccurlyeq_K \theta_Y$, this is equivalent to $T'(\bar{x}, v) \in -K$. Assume, by the way of contradiction, that $T'(\bar{x}, v) \notin -K$, this implies that $T'(\bar{x}, v) \in Y \setminus (-K)$. Notice that $-K$ is also a closed convex pointed cone in $Y$, which implies that $Y \setminus (-K)$ is an open subset of $Y$. From the assumption $T'(\bar{x}, v) \in Y \setminus (-K)$, there is $\varepsilon_1 > 0$ such that, $B_Y^o(T'(\bar{x}, v), \varepsilon_1) \subseteq Y \setminus (-K)$. This induces

$$B_Y^o(T'(\bar{x}, v), \varepsilon_1) \cap (-K) = \emptyset. \tag{2.10}$$

(2.10) implies that, for any $y \in Y$, we have

$$y \in -K \quad \Longrightarrow \quad \|y - T'(\bar{x}, v)\|_Y \geq \varepsilon_1. \tag{2.11}$$

By the Gâteaux directional differentiability of $T$ at point $\bar{x}$ along the given (fixed) direction $v$, by definition and by the condition that $A$ is an open subset in $X$ and $\bar{x} \in A$, for the above given $\varepsilon_1 > 0$, there is $\delta_1 > 0$ such that, for real $t$, we have

$$0 < |t| < \delta_1 \quad \Longrightarrow \bar{x} + tv \in A \ \text{ and } \ \left\|\frac{T(\bar{x}+tv)-T(\bar{x})}{t} - T'(\bar{x}, v)\right\|_Y < \varepsilon_1. \tag{2.12}$$

In particular, let $t < 0$ in (2.12), we get

$$0 < -t < \delta_1 \quad \Longrightarrow \quad \bar{x} - tv \in A \ \text{ and } \ \left\|\frac{T(\bar{x}-tv)-T(\bar{x})}{-t} - T'(\bar{x}, v)\right\|_Y < \varepsilon_1. \tag{2.13}$$

By the condition of this theorem that $T(\bar{x} - tv) - T(\bar{x}) \succcurlyeq_K \theta_Y$, that is, $T(\bar{x} - tv) - T(\bar{x}) \in K$, we must

have that

$$0 < -t < \delta_1 \quad \Longrightarrow \quad \frac{T(\bar{x}-tv)-T(\bar{x})}{-t} \in -K. \tag{2.14}$$

For $0 < -t < \delta_1$, in (2.14), (2.13) and (2.12), this creates a contradiction between (2.14), (2.13) and (2.11). Hence $T'(\bar{x}, v) \in -K$. This is

$$T'(\bar{x}, v) \preccurlyeq_K \theta_Y. \tag{2.15}$$

By (2.9) and (2.15) and by the fact that $\preccurlyeq_K$ is a partial order on $Y$, this proves $T'(\bar{x}, v) = \theta_Y$. □

Suppose that $T$ is Gâteaux differentiable at $\bar{x}$. More precisely speaking, in Corollary 2.5 below, by Theorem 2.4, we prove the following property.

$$\bar{x} \text{ is an } \preccurlyeq_K\text{-extreme point of } T \quad \Longrightarrow \quad \bar{x} \text{ is a generalized critical point of } T.$$

**Corollary 2.5**. *Let $\bar{x} \in A$. Suppose that $T$ is Gâteaux differentiable at $\bar{x}$. Then we have*

(i) *If $T$ has an (absolute) $\preccurlyeq_K$-maximum ((absolute) $\preccurlyeq_K$-minimum) value in $A$ at $\bar{x}$, then $\bar{x}$ is a generalized critical point of $T$ with*

$$T'(\bar{x})(v) = \theta_Y, \text{ for every } v \in X\backslash\{\theta_X\}.$$

*That is,* $$E(T) \subseteq dT^0.$$

(ii) *Conversely, if*

$$T'(\bar{x})(v) \neq \theta_Y, \text{ for some } v \in X\backslash\{\theta_X\},$$

*then $T$ takes neither $\preccurlyeq_K$-maximum, nor $\preccurlyeq_K$-minimum value at $\bar{x}$.*

*Proof.* This corollary follows from Theorem 2.4 and the proof is omitted here. □

**Corollary 2.6**. *Suppose that $T$ is Fréchet differentiable at $\bar{x}$. Then we have*

(i) *If $T$ has (absolute) $\preccurlyeq_K$-maximum ((absolute) $\preccurlyeq_K$-minimum) value in $A$ at $\bar{x}$, then*

$$\nabla T(\bar{x})(v) = T'(\bar{x})(v) = \theta_Y, \text{ for every } v \in X\backslash\{\theta_X\}.$$

(ii) *Inversely, if*

$$\nabla T(\bar{x})(v) \neq \theta_Y, \text{ for some } v \in X\backslash\{\theta_X\},$$

*then $T$ takes neither $\preccurlyeq_K$-maximum, nor $\preccurlyeq_K$-minimum value at point $\bar{x}$.*

*Proof.* By the fact that in Banach spaces, Fréchet differentiability of $T$ at $\bar{x}$ induces the Gâteaux differentiability of $T$ at $\bar{x}$, which satisfies $\nabla T(\bar{x}) = T'(\bar{x})$. Then, this corollary follows from Theorem 2.4 and Corollary 2.5 immediately. □

In particular, let $(\mathbb{R}, \leq)$ be the ordinary special totally ordered Hilbert space (in which the ordinary order $\leq$ is induced by $[0, \infty)$]. Let $(Y, \tau_Y, \preccurlyeq_K) = (\mathbb{R}, \leq)$ in Theorem 2.4, Corollary 2.5 and Corollary 2.6, we have the following ordinary maximum and minimum properties.

**Corollary 2.7**. *Let* $T: A \to \mathbb{R}$ *be a real valued function. Let* $\bar{x} \in A$. *We have that*

(i) *Let* $v \in X \backslash \{\theta_X\}$. *If* $T$ *has a maximum (minimum) value at* $\bar{x}$ *along direction* $v$ *and* $T$ *is Gâteaux directionally differentiable at* $\bar{x}$ *along direction* $v$, *then*

$$T'(\bar{x}, v) = 0.$$

(ii) *If* $T$ *has a maximum (minimum) value in* $A$ *at* $\bar{x}$ *and* $T$ *is Gâteaux differentiable at* $\bar{x}$, *then*

$$T'(\bar{x})(v) = 0, \text{ for every } v \in X \backslash \{\theta_X\}.$$

(iii) *Suppose that* $T$ *is Gâteaux differentiable at* $\bar{x}$. *If* $T$ *satisfies*

$$T'(\bar{x})(v) \neq 0, \text{ for some } v \in X \backslash \{\theta_X\},$$

*then* $T$ *takes neither maximum, nor minimum value at point* $\bar{x}$.

*Proof.* This corollary follows from Theorem 2.4 and the proof is omitted here. □

**Corollary 2.8**. *Let* $\bar{x} \in A$. *Let* $T: A \to \mathbb{R}$ *be a real valued functional. Suppose that* $T$ *is Fréchet differentiable at* $\bar{x}$. *Then*

(i) *If* $T$ *has a maximum (minimum) value in* $A$ *at point* $\bar{x}$, *then*

$$\nabla T(\bar{x})(v) = T'(\bar{x})(v) = 0, \text{ for every } v \in X \backslash \{\theta_X\}.$$

(ii) *Conversely, if*

$$\nabla T(\bar{x})(v) \neq 0, \text{ for some } v \in X \backslash \{\theta_X\},$$

*then* $T$ *takes neither maximum, nor minimum value at point* $\bar{x}$.

*Proof.* This corollary follows from Theorem 2.4 and the proof is omitted here. □

### 2.2. Order Monotone Property in Partially Ordered Banach spaces

In calculus, the monotone property of a differentiable real function on an open interval is described by the positive or negative of derivatives of the considered function. In [27], this property was extended to single-valued mappings in partially ordered Hausdorff topological vector spaces. In this subsection, we consider this property in partially ordered Banach spaces, which are special cases of partially ordered Hausdorff topological vector spaces and have some special properties.

In this subsection, unless otherwise stated, we always let $(X, \|\cdot\|_X, \preccurlyeq_C)$ and $(Y, \|\cdot\|_Y, \preccurlyeq_K)$ be partially ordered Banach spaces, in which the partial orders $\preccurlyeq_C$ and $\preccurlyeq_K$ are generated by nontrivial closed convex and pointed cones $C$ and $K$ in $X$ and $Y$, respectively, with $C \neq \{\theta_X\}$ and $K \neq \{\theta_Y\}$. Let $A$ be a nonempty convex open subset in $X$. Let $T: A \to Y$ be a single-valued mapping. We review the following definitions, which are used in [27]. If for any $x_1, x_2 \in A$,

$$x_1 \preccurlyeq_C x_2 \implies T(x_1) \preccurlyeq_K T(x_2),$$

then, $T$ is said to be $\preccurlyeq_C$-$\preccurlyeq_K$ increasing on $A$, or $T$ is said to be order increasing on $A$, if it is not confused. In particular, when $(Y, \tau_Y, \preccurlyeq_K) = (\mathbb{R}, \leq)$, then $T$ is said to be increasing on $A$. The concept of order decreasing can be similarly defined.

**Theorem 2.9**. *Suppose that $T$ is Gâteaux differentiable on $A$. If $T$ is order increasing on $A$, then for any $\bar{x} \in A$, we have*

$$T'(\bar{x})(v) \succcurlyeq_K \theta_Y, \text{ for any } v \in C \backslash \{\theta_X\}. \tag{2.16}$$

*In particular, if $T$ is Fréchet differentiable on $A$ and $T$ is order increasing on $A$, then for any $\bar{x} \in A$,*

$$\nabla T(\bar{x})(v) = T'(\bar{x})(v) \succcurlyeq_K \theta_Y, \text{ for any } v \in C \backslash \{\theta_X\}.$$

*Proof*. We first prove (2.16). Suppose that $T$ is order increasing on $A$. Since $A$ is a nonempty convex open subset in $X$, for the given $\bar{x} \in A$ and $v \in C \backslash \{\theta_X\}$, there is $\delta_0 > 0$ such that, for real number $t$, we have

$$|t| < \delta_0 \quad \Longrightarrow \quad \bar{x} + tv \in A.$$

By $v \in C \backslash \{\theta_X\}$, we have that

$$t > 0 \quad \Longrightarrow \quad \bar{x} + tv \succcurlyeq_C \bar{x} \quad \text{and} \quad t < 0 \quad \Longrightarrow \quad \bar{x} + tv \preccurlyeq_C \bar{x}.$$

By the condition that $T$ is order increasing on $A$, we always have

$$\frac{T(\bar{x}+tv)-T(\bar{x})}{t} \succcurlyeq_K \theta_Y, \text{ for any } t \text{ with } |t| < \delta_0. \tag{2.17}$$

Then, we will show that $T'(\bar{x})(v) \succcurlyeq_K \theta_Y$, this is equivalent to $T'(\bar{x})(v) \in K$. Assume, by the way of contradiction, that $T'(\bar{x})(v) \notin K$, this implies that $T'(\bar{x})(v) \in Y \backslash K$. Since $Y \backslash K$ is an open subset of $Y$, by the assumption that $T'(\bar{x})(v) \in Y \backslash K$, there is $\varepsilon > 0$ such that $B_Y^o(T'(\bar{x})(v), \varepsilon) \subseteq Y \backslash K$, which satisfies

$$B_Y^o(T'(\bar{x})(v), \varepsilon) = \{y \in Y : \|y - T'(\bar{x})(v)\|_Y < \varepsilon\}. \tag{2.18}$$

Meanwhile, the condition $B_Y^o(T'(\bar{x})(v), \varepsilon) \subseteq Y \backslash K$ induces

$$B_Y^o(T'(\bar{x})(v), \varepsilon) \cap K = \emptyset. \tag{2.19}$$

(2.18) and (2.19) together imply that, for any $y \in Y$, we have

$$y \in K \quad \Longrightarrow \quad \|y - T'(\bar{x})(v)\|_Y \geq \varepsilon. \tag{2.20}$$

Since $T$ is Gâteaux differentiable at point $\bar{x}$, then $T$ is Gâteaux directional differentiable at point $\bar{x}$ along the given (fixed) direction $v$. By definition and by the condition that $A$ is a convex open subset in $X$ and $\bar{x} \in A$, for the above given $\varepsilon > 0$, there is $\delta > 0$ such that, for real $t$, we have

$$0 < |t| < \delta \quad \Longrightarrow \bar{x} + tv \in A \text{ and } \left\|\frac{T(\bar{x}+tv)-T(\bar{x})}{t} - T'(\bar{x})(v)\right\|_Y < \varepsilon. \tag{2.22}$$

By (2.17), we have that

$$\frac{T(\bar{x}+tv)-T(\bar{x})}{t} \in K, \text{ for any } t \text{ with } 0 < |t| < \delta_0. \tag{2.23}$$

By (2.23), it creates a contradiction between (2.22) and (2.20). Hence, we must have

$$T'(\bar{x})(v) \succcurlyeq_K \theta_Y.$$

This proves (2.16). Then, by using (2.16), the second part of this theorem follows immediately from the connection between Fréchet differentiability and Gâteaux differentiability in Banach spaces. □

**Theorem 2.10**. *Let $T: A \to \mathbb{R}$ be a real valued function. Suppose that T is Gâteaux differentiable on A. If T is increasing on A, then for any $\bar{x} \in A$, we have*

$$T'(\bar{x})(v) \geq 0, \text{ for any } \in C\backslash\{\theta_X\}.$$

*In particular, if T is Fréchet differentiable on A and T is increasing on A, then for any $\bar{x} \in A$, we have*

$$\nabla T(\bar{x})(v) = T'(\bar{x})(v) \geq 0, \text{ for any } v \in C\backslash\{\theta_X\}.$$

*Proof.* This corollary follows from Theorem 2.9 immediately, in which, let $(Y, \tau_Y, \preccurlyeq_K) = (\mathbb{R}, \leq)$. □

In next section, we provide some examples to verify Theorem 2.9.

## 3. Gâteaux and Fréchet Derivatives in $l_p$ with $1 < p < \infty$

### 3.1. Some Explicit Formulas of Gâteaux and Fréchet Derivatives in $l_p$

Throughout this paper, unless otherwise stated, we always let $\mathbb{N}$, $\mathbb{R}$ and $\mathbb{R}_+$ respectively denote the set of nonnegative integers, the set of real numbers and the set of nonnegative real numbers. For given positive numbers $p$ and $q$ with $1 < p, q < \infty$ satisfying $\frac{1}{p} + \frac{1}{p} = 1$, let $(l_p, \|\cdot\|_p)$ be the common Banach space with topological dual space $(l_q, \|\cdot\|_q)$, which have the same origins $\theta_q = \theta_p = (0, 0, \ldots)$. Let $\langle \cdot, \cdot \rangle$ denote the pairing between $l_q$ and $l_p$. We define a mapping $\mathcal{P}: l_p \to \mathbb{R}_+$ as follows

$$\mathcal{P}(x) = \sup\{|t_n|: n = 1, 2, \ldots\}, \text{ for any } x = \{t_n\}_{n=1}^{\infty} \in l_p. \tag{3.1}$$

It is clear that $\mathcal{P}$ is well-defined on $l_p$ and it satisfies that,

$$0 \leq \mathcal{P}(x) < \infty, \text{ for any } x = \{t_n\}_{n=1}^{\infty} \in l_p.$$

Furthermore, by the summability of the $p^{th}$ power, for any $x = \{t_n\}_{n=1}^{\infty} \in l_p$, we have that,

$$\mathcal{P}(x) = \sup\{|t_n|: n = 1, 2, \ldots\} = \max\{|t_n|: n = 1, 2, \ldots\}.$$

Let $m$ be a positive integer. We define the $m^{th}$ power operator $Q^m: l_p \to l_p$ as follows.

$$Q^m(x) = \{t_n^m\}_{n=1}^{\infty}, \text{ for any } x = \{t_n\}_{n=1}^{\infty} \in l_p. \tag{3.2}$$

One could easily check that the $m^{th}$ power operator $Q^m$ is well-defined from $l_p$ to $l_p$.

**Theorem 3.1**. *Let m be a positive integer. Then, the $m^{th}$ power operator $Q^m$ is Fréchet differentiable on $l_p$ such that for given $\bar{x} = (\bar{t}_1, \bar{t}_2, \ldots) \in l_p$, the Fréchet derivative $\nabla(Q^m)(\bar{x}): l_p \to l_p$ of $Q^m$ at $\bar{x}$ satisfies the following properties.*

(i) *If $m > 1$, then*

$$\nabla(Q^m)(\bar{x}) = \begin{pmatrix} m\bar{t}_1^{m-1} & 0 & \dots \\ 0 & m\bar{t}_2^{m-1} & 0 \\ \vdots & 0 & \ddots \end{pmatrix}. \tag{3.3}$$

*Here, $\nabla(Q^m)(\bar{x})$ is a pointwise multiplication operator on $l_p$ such that, for $x = \{t_n\}_{n=1}^{\infty} \in l_p$,*

$$\nabla(Q^m)(\bar{x})(x) = \{t_n\}_{n=1}^{\infty} \begin{pmatrix} m\bar{t}_1^{m-1} & 0 & \dots \\ 0 & m\bar{t}_2^{m-1} & 0 \\ \vdots & 0 & \ddots \end{pmatrix} = \{m\bar{t}_n^{m-1} t_n\}_{n=1}^{\infty} \in l_p. \tag{3.4}$$

(ii) *If $m = 1$, then for given $\bar{x} = (\bar{t}_1, \bar{t}_2, \dots) \in l_p$, the Fréchet derivative $\nabla(Q)(\bar{x})$ is the following identity mapping on $l_p$*

$$\nabla(Q)(\bar{x}) = \begin{pmatrix} 1 & 0 & \dots \\ 0 & 1 & 0 \\ \vdots & 0 & \ddots \end{pmatrix}. \tag{3.5}$$

*This satisfies that*

$$\nabla(Q)(\bar{x})(x) = x = \{t_n\}_{n=1}^{\infty}, \text{ for every } x = \{t_n\}_{n=1}^{\infty} \in l_p. \tag{3.6}$$

*Proof.* For an arbitrarily given $\bar{x} = (\bar{t}_1, \bar{t}_2, \dots) \in l_p$, let $\nabla(Q^m)(\bar{x})$ be represented by (3.3), which acts as a pointwise multiplication operator on $l_p$ as given in (3.4). Then, for any $u = \{s_n\}_{n=1}^{\infty} \in l_p \backslash \{\theta_p\}$ with $0 < \|u\|_p < 1$, we calculate

$$\frac{\|Q^m(\bar{x}+u) - Q^m(\bar{x}) - \nabla(Q^m)(\bar{x})(u)\|_p}{\|u\|_p}$$

$$= \frac{\left\|\{(\bar{t}_n+s_n)^m - \bar{t}_n^m - m\bar{t}_n^{m-1}s_n\}_{n=1}^{\infty}\right\|_p}{\|u\|_p}$$

$$= \frac{\left\|\{\binom{m}{2}\bar{t}_n^{m-2}s_n^2 + \binom{m}{3}\bar{t}_n^{m-3}s_n^3 + \dots + \binom{m}{m}s_n^m\}_{n=1}^{\infty}\right\|_p}{\|u\|_p}$$

$$\leq \frac{\left(\sum_{n=1}^{\infty}\left(\binom{m}{2}|\bar{t}_n|^{m-2}|s_n|^2 + \binom{m}{3}|\bar{t}_n|^{m-3}|s_n|^3 + \dots + \binom{m}{m}|s_n|^m\right)^p\right)^{\frac{1}{p}}}{\|u\|_p}$$

$$\leq \frac{\left(\sum_{n=1}^{\infty}\left(\binom{m}{2}\left(\mathcal{P}(\bar{x})\right)^{m-2}|s_n|^2 + \binom{m}{3}\left(\mathcal{P}(\bar{x})\right)^{m-3}|s_n|^3 + \dots + \binom{m}{m}|s_n|^m\right)^p\right)^{\frac{1}{p}}}{\|u\|_p}$$

$$\leq \frac{\left(\sum_{n=1}^{\infty}\left(|s_n|^2\left(\binom{m}{2}\left(\mathcal{P}(\bar{x})\right)^{m-2} + \binom{m}{3}\left(\mathcal{P}(\bar{x})\right)^{m-3} + \dots + \binom{m}{m}\right)\right)^p\right)^{\frac{1}{p}}}{\|u\|_p}$$

$$\leq \frac{\left(\sum_{n=1}^{\infty}\left(|s_n|^2\left(1+\mathcal{P}(\bar{x})\right)^m\right)^p\right)^{\frac{1}{p}}}{\|u\|_p}$$

$$= \frac{\left(1+\mathcal{P}(\bar{x})\right)^m \left(\sum_{n=1}^{\infty}\left(|s_n|^2\right)^p\right)^{\frac{1}{p}}}{\|u\|_p}$$

$$= \frac{\left(1+\mathcal{P}(\bar{x})\right)^m \left(\sum_{n=1}^{\infty}|s_n|^p|s_n|^p\right)^{\frac{1}{p}}}{\|u\|_p}$$

$$\leq \frac{\left(1+\mathcal{P}(\bar{x})\right)^m \left(\sum_{n=1}^{\infty}|s_n|^p\left(\sum_{i=1}^{\infty}|s_i|^p\right)\right)^{\frac{1}{p}}}{\|u\|_p}$$

$$= \frac{\left(1+\mathcal{P}(\bar{x})\right)^m \left(\sum_{n=1}^{\infty}|s_n|^p\right)^{\frac{1}{p}}\left(\left(\sum_{i=1}^{\infty}|s_i|^p\right)\right)^{\frac{1}{p}}}{\|u\|_p}$$

$$= \left(1+\mathcal{P}(\bar{x})\right)^m \|u\|_p.$$

This implies that

$$\lim_{u \to \theta_X} \frac{\|Q^m(\bar{x}+u) - Q^m(\bar{x}) - \nabla(Q^m)(\bar{x})(u)\|_q}{\|u\|_p} = 0. \qquad \square$$

Let $m$ be a positive integer. Let $a_m$, $a_{m-1}$, … $a_1$ be real numbers. Define a polynomial type operator $Q_m: l_p \to l_p$, for $x = \{t_n\}_{n=1}^{\infty} \in l_p$, by

$$Q_m(x) = a_m x^m + a_{m-1} x^{m-1} + \ldots + a_1 x = \left\{\sum_{i=1}^{m} a_i t_n^i\right\}_{n=1}^{\infty}. \tag{3.7}$$

**Theorem 3.2**. *Let m be a positive integer and let $Q_m$ be a polynomial type operator on $l_p$ defined by* (3.7). *Then, $Q_m$ is Fréchet differentiable on $l_p$ such that for any given $\bar{x} = \{\bar{t}_n\}_{n=1}^{\infty} \in l_p$, the Fréchet derivative $\nabla Q_m(\bar{x}): l_p \to l_p$ satisfies the following properties.*

(i) *If m > 1, we have*

$$\nabla Q_m(\bar{x}) = \begin{pmatrix} \sum_{i=1}^{m} i a_i \bar{t}_1^{i-1} & 0 & \cdots \\ 0 & \sum_{i=1}^{m} i a_i \bar{t}_2^{i-1} & 0 \\ \vdots & 0 & \ddots \end{pmatrix}. \tag{3.8}$$

*Here, $\nabla Q_m(\bar{x})$ satisfies that, for every $x = \{t_n\}_{n=1}^{\infty} \in l_p$,*

$$\nabla Q_m(\bar{x})(x) = \{t_n\}_{n=1}^{\infty} \begin{pmatrix} \sum_{i=1}^{m} i a_i \bar{t}_1^{i-1} & 0 & \cdots \\ 0 & \sum_{i=1}^{m} i a_i \bar{t}_2^{i-1} & 0 \\ \vdots & 0 & \ddots \end{pmatrix} = \left\{\sum_{i=1}^{m} i a_i \bar{t}_n^{i-1} t_n\right\}_{n=1}^{\infty} \in l_p.$$

(ii) *When m = 1, we have*

$$\nabla Q_1(\bar{x}) = \begin{pmatrix} a_1 & 0 & \cdots \\ 0 & a_1 & 0 \\ \vdots & 0 & \ddots \end{pmatrix}.$$

*For every $x = \{t_n\}_{n=1}^{\infty} \in l_p$,*

$$\nabla Q_1(\bar{x}) x = a_1 x = \{a_1 t_n\}_{n=1}^{\infty} \in l_p.$$

*Proof*. This theorem follows from Theorem 3.1 and the linearity of Fréchet derivatives immediately. □

**Corollary 3.3**. *Let m be a positive integer and let $Q_m$ be a polynomial type operator on $l_p$ defined by* (3.7). *Then, $Q_m$ is Gâteaux differentiable on $l_p$ such that, for any $\bar{x} = (\bar{t}_1, \bar{t}_2, \dots) \in l_p$, the Gâteaux derivative $(Q_m)'(\bar{x})$ satisfies the following properties.*

(i) *If m > 1, we have*

$$(Q_m)'(\bar{x}) = \nabla Q_m(\bar{x}) = \begin{pmatrix} \sum_{i=1}^{m} i a_i \bar{t}_1^{i-1} & 0 & \dots \\ 0 & \sum_{i=1}^{m} i a_i \bar{t}_2^{i-1} & 0 \\ \vdots & 0 & \ddots \end{pmatrix}. \tag{3.9}$$

*Here, $(Q_m)'(\bar{x})$ is represented by an $\infty \times \infty$ matrix that defines a pointwise multiplication operator on $l_p$ satisfying that, for every $v = \{t_n\}_{n=1}^{\infty} \in l_p \backslash \{\theta_p\}$,*

$$(Q_m)'(\bar{x})(v) = \nabla Q_m(\bar{x})(v) = \{t_n\}_{n=1}^{\infty} \begin{pmatrix} \sum_{i=1}^{m} i a_i \bar{t}_1^{i-1} & 0 & \dots \\ 0 & \sum_{i=1}^{m} i a_i \bar{t}_2^{i-1} & 0 \\ \vdots & 0 & \ddots \end{pmatrix} = \left\{\sum_{i=1}^{m} i a_i \bar{t}_n^{i-1} t_n\right\}_{n=1}^{\infty}.$$

(ii) *When m = 1, we have*

$$(Q_1)'(\bar{x}) = \nabla Q_1(\bar{x}) = \begin{pmatrix} a_1 & 0 & \dots \\ 0 & a_1 & 0 \\ \vdots & 0 & \ddots \end{pmatrix}.$$

*For every $v = \{t_n\}_{n=1}^{\infty} \in l_p \backslash \{\theta_p\}$, we have*

$$(Q_m)'(\bar{x})(v) = \nabla Q_1(\bar{x}) v = a_1 v = \{a_1 t_n\}_{n=1}^{\infty} \in l_p.$$

Recall that $d(Q_m)^0$ is the set of generalized critical point of $Q_m$.

**Corollary 3.4**. *Let m be a positive integer and let $Q_m$ be a polynomial type operator on $l_p$ defined by* (3.7). *Then, $Q_m$ has the following properties.*

(i) *If m > 1 with $a_m \neq 0$,* and $a_1 = 0$, *then*

$$\nabla(Q_m)^0 = d(Q_m)^0$$
$$= \{\bar{x} = (\bar{t}_1, \bar{t}_2, \dots) \in l_p : \textstyle\sum_{i=2}^{m} i a_i \bar{t}_n^{i-1} = 0, \text{for } n = 1, 2, \dots\}.$$

*More precisely, define the real solution set A = $\{s \in \mathbb{R} : \sum_{i=2}^{m} i a_i s^{i-1} = 0\}$. If $\{0\} \subsetneq A$, then $\{\theta_p\} \subsetneq \nabla(Q_m)^0$; and*

$$\bar{x} = (\bar{t}_1, \bar{t}_2, \dots) \in \nabla(Q_m)^0 \backslash \{\theta_p\}$$

*if and only if, there is a nonempty finite subset $\mathcal{F} \subsetneq \mathbb{N}$ such that*

$$\bar{t}_n = \begin{cases} s, & \text{if } n \in \mathcal{F}, \text{for some } s \in A \backslash \{0\}, \\ 0, & \text{if } n \in \mathbb{N} \backslash \mathcal{F}, \end{cases} \quad \text{for } n = 1, 2, \dots .$$

(ii) *If m ≥ 1 with $a_1 \neq 0$, then*

$$\nabla(Q_m)^0 = d(Q_m)^0 = \emptyset.$$

*Proof.* Proof of (i). By Corollary 3.3, for $\bar{x} = (\bar{t}_1, \bar{t}_2, \dots) \in l_p$, by $a_m \neq 0$ and $a_1 = 0$, we have

$$\nabla Q_m(\bar{x}) = \begin{pmatrix} \sum_{i=2}^{m} i a_i \bar{t}_1^{i-1} + a_1 & 0 & \dots \\ 0 & \sum_{i=2}^{m} i a_i \bar{t}_2^{i-1} + a_1 & 0 \\ \vdots & 0 & \ddots \end{pmatrix}$$

$$= \begin{pmatrix} \sum_{i=2}^{m} i a_i \bar{t}_1^{i-1} & 0 & \dots \\ 0 & \sum_{i=2}^{m} i a_i \bar{t}_2^{i-1} & 0 \\ \vdots & 0 & \ddots \end{pmatrix}$$

By definition, for $\bar{x} = (\bar{t}_1, \bar{t}_2, \dots) \in l_p$, we have that

$$\bar{x} \in \nabla(Q_m)^0 \Longrightarrow \nabla Q_m(\bar{x}) = \theta_{pp}.$$

Here, $\theta_{pp}: l_p \to l_p$ denotes the zero continuous and linear operator from $l_p$ to $l_p$, which implies that,

$$\nabla(Q_m)^0 = \{\bar{x} = (\bar{t}_1, \bar{t}_2, \dots) \in l_p : \textstyle\sum_{i=2}^{m} i a_i \bar{t}_n^{i-1} = 0, \text{for } n = 1, 2, \dots\}.$$

Proof of (ii). There are two cases for consideration.

Case 1. Since $m = 1$ with $a_1 \neq 0$. In this case, by part (ii) in Corollary 3.3, for every $v = \{t_n\}_{n=1}^{\infty} \in l_p \backslash \{\theta_p\}$,

$$(Q_m)'(\bar{x})(v) = \nabla Q_1(\bar{x})v = a_1 v \in l_p.$$

This proves part (ii) immediately under conditions $m = 1$ and $a_1 \neq 0$.

Case 2. $m > 1$ with $a_1 \neq 0$. By Corollary 3.3, for $\bar{x} = (\bar{t}_1, \bar{t}_2, \dots) \in l_p$, by $a_1 \neq 0$, we have

$$\nabla Q_m(\bar{x}) = \begin{pmatrix} \sum_{i=2}^{m} i a_i \bar{t}_1^{i-1} + a_1 & 0 & \dots \\ 0 & \sum_{i=2}^{m} i a_i \bar{t}_2^{i-1} + a_1 & 0 \\ \vdots & 0 & \ddots \end{pmatrix}$$

Similarly to the proof of part (i), under the condition of (ii), we have

$$\nabla(Q_m)^0 = \{\bar{x} = (\bar{t}_1, \bar{t}_2, \dots) \in l_p : \textstyle\sum_{i=2}^{m} i a_i \bar{t}_n^{i-1} + a_1 = 0, \text{for } n = 1, 2, \dots\} = \emptyset.$$

The above equation follows from the properties that the system $\sum_{i=2}^{m} i a_i \bar{t}_n^{i-1} + a_1 = 0$, for $n = 1, 2, \dots$ has at most $m$ nonzero solutions. By $(\bar{t}_1, \bar{t}_2, \dots) \in l_p$, it is impossible for any $\bar{x} \in l_p$ to satisfy this system. This proves case (ii). □

### 3.2. Ordered Optimization in Partially Ordered Banach Space $l_p$

Let $(l_p, \|\cdot\|_p)$, for $1 < p < \infty$ with origin $\theta_p$, be the Banach space studied in subsection 3.1. Let $K$ be the positive cone in $(l_p, \|\cdot\|_p)$ that is defined by

$$K = \left\{u = \{s_n\}_{n=1}^{\infty} \in l_p : s_n \geq 0, \text{for each } n = 1, 2, \dots\right\}.$$

It is clear that $K$ is a nonempty closed convex and pointed cone in $l_p$ with $K \neq \{\theta_p\}$. Let $\preccurlyeq_K$ be the partial order on $l_p$ induced by $K$, by which $(l_p, \|\cdot\|_p, \preccurlyeq_K)$ is a partially ordered Banach space. Then, by Corollary 3.4 and part (ii) in Corollary 2.6, we investigate the $\preccurlyeq_K$-extrema of the polynomial type operator $Q_m$on $l_p$.

**Corollary 3.5**. *Let m be a positive integer and let $Q_m$ be the polynomial type operator on $l_p$ defined by* (3.7). *Then, $Q_m$ has the following properties.,*

(i) *If $m > 1$ with $a_m \neq 0$ and $a_1 = 0$, then $\theta_p$ is one of the possible $\preccurlyeq_K$-extrema point of $Q_m$;*
(ii) *If $m \geq 1$ with $a_1 \neq 0$, then $Q_m$ does not have $\preccurlyeq_K$-extrema.*

*Proof.* This corollary follows from Corollary 3.4 and part (ii) in Corollary 2.6 immediately. □

We first give an example in $(l_p, \|\cdot\|_p, \preccurlyeq_K)$ below to verify Theorem 2.9. That is, if a given mapping $T$ is Gâteaux differentiable on a convex and open subset $A$, then

$$T \text{ is } \preccurlyeq_K\text{-increasing on } A \Longrightarrow T'(\bar{x})(v) \succcurlyeq_K \theta_Y, \text{ for any } \bar{x} \in A \text{ and for any } v \in K\backslash\{\theta_X\}. \tag{3.10}$$

**Example 3.6**. Let $Q^3$ be the 3rd power operator on $l_p$ defined by

$$Q^3(u) = \{s_n^3\}_{n=1}^{\infty} \in l_p, \text{ for any } u = \{s_n\}_{n=1}^{\infty} \in l_p.$$

By (3.5), for given $\bar{x} = (\bar{t}_1, \bar{t}_2, \ldots) \in l_p$, the Fréchet derivative $\nabla(Q^3)(\bar{x})\colon l_p \to l_p$ is given by

$$\nabla(Q^3)(\bar{x}) = \begin{pmatrix} 3\bar{t}_1^2 & 0 & \ldots \\ 0 & 3\bar{t}_2^2 & 0 \\ \vdots & 0 & \ddots \end{pmatrix}. \tag{3.11}$$

For any $x = \{t_n\}_{n=1}^{\infty} \in l_p$, by (3.6) and (3.11), we have

$$\nabla(Q^3)(\bar{x})(x) = \{t_n\}_{n=1}^{\infty} \begin{pmatrix} 3\bar{t}_1^2 & 0 & \ldots \\ 0 & 3\bar{t}_2^2 & 0 \\ \vdots & 0 & \ddots \end{pmatrix} = \{3\bar{t}_n^2 t_n\}_{n=1}^{\infty} \in l_p. \tag{3.11}$$

For any $x = \{t_n\}_{n=1}^{\infty}$ and $y = \{s_n\}_{n=1}^{\infty} \in l_p$, by definition, we have that

$$x \preccurlyeq_K y \text{ if and only if } t_n \leq s_n, \text{ for all } n = 1, 2, \ldots . \tag{3.12}$$

Notice that, for each $n = 1, 2, \ldots$, $s_n^3 - t_n^3 = (s_n - t_n)(s_n^2 + s_n t_n + t_n^2)$ and $(s_n^2 + s_n t_n + t_n^2) \geq 0$. By (3.12), this implies that, for any $x, y \in l_p$,

$$x \preccurlyeq_K y \text{ if and only if } Q^3(x) \preccurlyeq_K Q^3(y). \tag{3.13}$$

Hence, $Q^3$ is an order increasing mapping on $l_p$. By (3.11) and (3.13), for $\bar{x} = (\bar{t}_1, \bar{t}_2, \ldots) \in l_p$, we have

$$\nabla Q^3(\bar{x})(v) = (Q^3)'(\bar{x})(v) = \{3\bar{t}_n^2 t_n\}_{n=1}^{\infty} \succcurlyeq_K \theta_Y, \text{ for any } v = \{t_n\}_{n=1}^{\infty} \in K\backslash\{\theta_X\},$$

This proves (3.10) and it verifies Theorem 2.9 for this 3rd power operator $Q^3$.

Next, we define a polynomial type power operator $Q_3$ on $l_p$, which satisfies the following properties,

For any $\bar{x} \in l_p$ and for any $v \in K\backslash\{\theta_p\}$, $(Q_3)'(\bar{x})(v) \succcurlyeq_K \theta_Y$ and $Q_3$ is $\preccurlyeq_K$-increasing on $l_p$. (4.14)

**Example 3.7**. Let $Q_3$ be a polynomial type operator on $l_p$ defined by

$$Q_3(u) = \frac{1}{3}u^3 + u^2 + u = \left\{\frac{1}{3}s_n^3 - s_n^2 + s_n\right\}_{n=1}^{\infty} \in l_p, \text{ for any } u = \{s_n\}_{n=1}^{\infty} \in l_p.$$

By (3.5), for given $\bar{x} = (\bar{t}_1, \bar{t}_2, \dots) \in l_p$, the Fréchet derivative $\nabla(Q_3)(\bar{x}): l_p \to l_p$ is satisfies that, for any $v = \{t_n\}_{n=1}^{\infty} \in K$, we have

$$\nabla(Q_3)(\bar{x})(v) = \{t_n\}_{n=1}^{\infty} \begin{pmatrix} (\bar{t}_1+1)^2 & 0 & \dots \\ 0 & (\bar{t}_2+1)^2 & 0 \\ \vdots & 0 & \ddots \end{pmatrix} = \{(\bar{t}_n+1)^2 t_n\}_{n=1}^{\infty} \in K.$$

This implies that

$$\nabla(Q_3)(\bar{x})(v) \succcurlyeq_K \theta_Y, \text{ for any } \bar{x} \in l_p \text{ and for any } v \in K\backslash\{\theta_X\}.$$

Let $x = \{t_n\}_{n=1}^{\infty}$ and $y = \{s_n\}_{n=1}^{\infty} \in l_p$. One can check that (3.14) is satisfied for this mapping.

$$x \preccurlyeq_K y \implies Q_3(x) \preccurlyeq_K Q_3(y).$$

The next example verifies Corollary 2.5. That is, if a given mapping is Gâteaux differentiable at $\bar{x}$, then

$\bar{x}$ is an $\preccurlyeq_K$-extreme point of given mapping $\implies$ $\bar{x}$ is a generalized critical point of this mapping. (3.15)

**Example 3.8**. For an arbitrarily given positive integer $m$, let $Q^{2m}$ be the $2m^{\text{th}}$ power operator on $l_p$ studied in section 3. $Q^{2m}$ is defined by

$$Q^{2m}(u) = \{s_n^{2m}\}_{n=1}^{\infty} \in l_p, \text{ for any } u = \{s_n\}_{n=1}^{\infty} \in l_p.$$

By definition, we have that

$$\{s_n^{2m}\}_{n=1}^{\infty} \in K, \text{ for any } u = \{s_n\}_{n=1}^{\infty} \in l_p.$$

This implies that

$$Q^{2m}(u) \succcurlyeq_K \theta_p, \text{ for any } u = \{s_n\}_{n=1}^{\infty} \in l_p.$$

Hence, this operator $Q^{2m}$ takes absolute $\preccurlyeq_K$-minimum value at $\theta_p$ in $l_p$. On the other hand, by Theorem 3.1, we have

$$\nabla(Q^{2m})(\theta_p) = \begin{pmatrix} 0 & 0 & \dots \\ 0 & 0 & 0 \\ \vdots & 0 & \ddots \end{pmatrix}.$$

This satisfies that

$$\nabla(Q^{2m})(\theta_p)(x) = \{t_n\}_{n=1}^{\infty} \begin{pmatrix} 0 & 0 & \dots \\ 0 & 0 & 0 \\ \vdots & 0 & \ddots \end{pmatrix} = \theta_p, \text{ for any } x = \{t_n\}_{n=1}^{\infty} \in l_p.$$

And, the point $\theta_p$ is a generalized critical point of $R^{2m}$ that verifies (3.15) and Corollary 2.4.

Next, by using the Banach space $l_p$, we will provide a counterexample to show that

$$\bar{x} \text{ is a generalized critical point of a given mapping } \nRightarrow \bar{x} \text{ is an ordered extreme of this mapping.} \quad (3.16)$$

**Example 3.9**. Consider that partially ordered Banach space $(l_p, \|\cdot\|_p, \preccurlyeq_K)$ with $1 < p < \infty$ defined in Example 3.8. Similar to Example 3.6, the 3rd power operator $Q^3$ on $l_p$ satisfies that

$$\nabla(Q^3)(\theta_p)(x) = \{t_n\}_{n=1}^{\infty} \begin{pmatrix} 0 & 0 & \dots \\ 0 & 0 & 0 \\ \vdots & 0 & \ddots \end{pmatrix} = \theta_p, \text{ for any } \ x = \{t_n\}_{n=1}^{\infty} \in l_p.$$

Hence, the point $\theta_p$ is a generalized critical point of $Q^3$. Let $w = \left\{\frac{1}{2^n}\right\}_{n=1}^{\infty} \in l_p \backslash \{\theta_p\}$. It is clear that $w \in K \backslash \{\theta_p\}$, that is, $w \succ_K \theta_p$. For this given $w \in K \backslash \{\theta_p\}$, we have that

$$t > 0 \implies Q^3(tw) = \left\{\frac{t^3}{2^{3n}}\right\}_{n=1}^{\infty} \in K \backslash \{\theta_p\}, \text{ that is } Q^3(tw) \succ_K \theta_p;$$

and

$$t < 0 \implies Q^3(tw) = \left\{\frac{t^3}{2^{3n}}\right\}_{n=1}^{\infty} \in -K \backslash \{\theta_p\}, \text{ that is } Q^3(tw) \prec_K \theta_p.$$

Hence, the origin $\theta_p$ is not an $\preccurlyeq_K$-extreme of the 3rd power operator $Q^3$ on $l_p$. This verifies (3.16). □

We provide the following example to demonstrate Corollaries 3.3 and 3.4.

**Example 3.10**. Let the partially ordered Banach space $(l_p, \|\cdot\|_p, \preccurlyeq_K)$ with $1 < p < \infty$ be defined in Example 3.8. Similar to Example 3.7, let $Q_4$ be a polynomial type operator on $l_p$ defined by

$$Q_4(u) = \frac{1}{4}u^4 - u^3 - 2u^2 = \left\{\frac{1}{4}s_n^4 - s_n^3 - 2s_n^2\right\}_{n=1}^{\infty} \in l_p, \text{ for any } u = \{s_n\}_{n=1}^{\infty} \in l_p.$$

By (3.5), for given $\bar{x} = (\bar{t}_1, \bar{t}_2, \dots) \in l_p$, the Fréchet derivative $\nabla(Q_4)(\bar{x}): l_p \to l_p$ is satisfies that

$$\nabla(Q_4)(\bar{x}) = \begin{pmatrix} \bar{t}_1^3 - 3\bar{t}_1^2 - 4\bar{t}_1 & 0 & \dots \\ 0 & \bar{t}_2^3 - 3\bar{t}_2^2 - 4\bar{t}_2 & 0 \\ \vdots & 0 & \ddots \end{pmatrix}.$$

By following Corollary 3.4, the real solution set for $\nabla(Q_4)$ is $A = \{0, -1, 4\}$. Then $\nabla(Q_4)^0$ is an infinite subset of $l_p$. We pick up some examples in $\nabla(Q_4)^0 \backslash \{\theta_p\}$. Let $\mathcal{F} = \{1, 5, 6, 8\}$. According to this $\mathcal{F}$, we construct 4 examples in $\nabla(Q_4)^0 \backslash \{\theta_p\}$ as follows.

$$(-1, 0, 0\ , 0, -1, 4, 0, -1, 0, 0, \dots); (4, 0, 0\ , 0, -1, 4, 0, 4, 0, 0, \dots);$$

$$(-1, 0, 0\ , 0, -1, -1, 0, -1, 0, 0, \dots); (4, 0, 0\ , 0, 4, 4, 0, 4, 0, 0, \dots).$$

## 4. Gâteaux and Fréchet Derivatives in *C*[0, 1]

### 4.1. Some Explicit Formulas of Gâteaux and Fréchet Derivatives in *C*[0, 1]

Let $(C[0, 1], \|\cdot\|_C)$ be the Banach space of all real continuous functions defined on [0, 1] equipped with the maximum norm $\|\cdot\|_C$ defined, for any $f \in C[0, 1]$, by

$$\|f\|_C = \max\{|f(x)|: x \in [0, 1]\}.$$

Let $\theta_C$ denote the null element of $C[0, 1]$. In this subsection, we investigate the explicit formulas of the Gâteaux and Fréchet derivatives of some single-valued mappings in $C[0, 1]$.

**Theorem 4.1**. *Let m be a positive integer and let* $P^m: C[0, 1] \to C[0, 1]$ *be the* $m^{th}$ *power operator on* $C[0, 1]$ *defined by*

$$P^m(f) = f^m, \text{ for each } f \in C[0, 1].$$

*Then,* $P^m$ *is Fréchet differentiable on* $C[0, 1]$ *such that for any given* $\bar{f} \in C[0, 1]$, *the Fréchet derivative* $\nabla P^m(\bar{f}): C[0, 1] \to C[0, 1]$ *satisfies that*

$$\nabla P^m(\bar{f})(g) = m\bar{f}^{m-1}g, \text{ for any } g \in C[0, 1]. \tag{4.1}$$

*Proof.* It is clear that (4.1) holds for $m = 1$. Hence, we may assume that $m > 1$. Let $\bar{f} \in C[0, 1]$ be arbitrarily given. Then, we calculate

$$\lim_{u \xrightarrow{C[0,1]} \theta_C} \left\| \frac{P^m(\bar{f}+u) - P^m(\bar{f}) - m\bar{f}^{m-1}u}{\|u\|_C} \right\|_C$$

$$= \lim_{u \xrightarrow{C[0,1]} \theta_C} \frac{\|(\bar{f}+u)^m - \bar{f}^m - m\bar{f}^{m-1}u\|_C}{\|u\|_C}$$

$$= \lim_{u \xrightarrow{C[0,1]} \theta_C} \frac{\left\|\binom{m}{2}\bar{f}^{m-2}u^2 + \binom{m}{3}\bar{f}^{m-3}u^3 + \binom{m}{m}u^m\right\|_C}{\|u\|_C}$$

$$\leq \lim_{u \xrightarrow{C[0,1]} \theta_C} \|u\|_C \left\|\binom{m}{2}\bar{f}^{m-2} + \binom{m}{3}\bar{f}^{m-3}u + \binom{m}{m}u^{m-2}\right\|_C$$

$$= 0.$$ □

Let $m$ be a positive integer and let $a_m, a_{m-1}, \ldots a_1$ be real numbers. Define a polynomial type operator $P_m: C[0, 1] \to C[0, 1]$, for $f \in C[0, 1]$, by

$$P_m(f) = a_m f^m + a_{m-1} f^{m-1} + \ldots + a_1 f. \tag{4.2}$$

**Corollary 4.2**. *Let m be a positive integer and let* $P_m$ *be the polynomial type operator on* $C[0, 1]$ *defined by* (4.2). *Then,* $P_m$ *has the following properties.*

(i) $P_m$ *is Fréchet differentiable on* $C[0, 1]$ *such that for any given* $\bar{f} \in C[0, 1]$, *the Fréchet derivative* $\nabla P_m(\bar{f}): C[0, 1] \to C[0, 1]$ *satisfies that*

$$\nabla P_m(\bar{f})(g) = \sum_{i=1}^{m} i a_i \bar{f}^{i-1} g, \text{ for any } g \in C[0, 1]. \tag{4.3}$$

(ii) $P_m$ *is Gâteaux differentiable on* $C[0, 1]$ *such that for any given* $\bar{f} \in C[0, 1]$, *the Gâteaux derivative* $(P_m)'(\bar{f})$ *satisfies that*

$$(P_m)'(\bar{f})(g) = \nabla P_m(\bar{f})(g) = \sum_{i=1}^{m} i a_i \bar{f}^{i-1} g, \textit{ for any } g \in C[0,1]. \tag{4.4}$$

*Proof.* By (4.1) in Theorem 4.1, the proof of this corollary is straight forward and it is omitted here. □

**Theorem 4.3**. *Define a single-valued mapping* $S: C[0,1] \to C[0,1]$ *by*

$$S(f) = \sin(f), \textit{ for each } f \in C[0,1].$$

*Then,* $S$ *is Fréchet differentiable on* $C[0, 1]$ *such that for any given* $\bar{f} \in C[0,1]$, *the Fréchet derivative* $\nabla S(\bar{f}): C[0,1] \to C[0,1]$ *satisfies that*

$$\nabla S(\bar{f})(g) = \mathrm{con}(\bar{f})g, \textit{ for any } g \in C[0,1].$$

*And therefore,* $S$ *is Gâteaux differentiable on* $C[0, 1]$ *and the Gâteaux derivative* $S'(\bar{f})$ *satisfies that*

$$S'(\bar{f})(g) = \nabla S(\bar{f})(g) = \mathrm{con}(\bar{f})g, \textit{ for any } g \in C[0,1]. \tag{4.5}$$

*Proof.* This theorem can be proved by using the chain rule of Fréchet derivatives. However, we give a direct proof of this theorem here. Let $\bar{f} \in C[0,1]$ be arbitrarily given. Then, we calculate

$$\begin{aligned}
&\lim_{u \xrightarrow{C[0,1]} \theta_C} \left\| \frac{\sin(\bar{f}+u) - \sin(\bar{f}) - \mathrm{con}(\bar{f})u}{\|u\|_C} \right\|_C \\
&= \lim_{u \xrightarrow{C[0,1]} \theta_C} \left\| \frac{\sin(\bar{f})\mathrm{con}(u) + \mathrm{con}(\bar{f})\sin(u) - \sin(\bar{f}) - \mathrm{con}(\bar{f})u}{\|u\|_C} \right\|_C \\
&\le \lim_{u \xrightarrow{C[0,1]} \theta_C} \left( \left\| \frac{\sin(\bar{f})\mathrm{con}(u) - \sin(\bar{f})}{\|u\|_C} \right\|_C + \left\| \frac{\mathrm{con}(\bar{f})\sin(u) - \mathrm{con}(\bar{f})u}{\|u\|_C} \right\|_C \right) \\
&\le \lim_{u \xrightarrow{C[0,1]} \theta_C} \left( \left\| \frac{\sin(\bar{f})\mathrm{con}(u) - \sin(\bar{f})}{u} \right\|_C + \left\| \frac{\mathrm{con}(\bar{f})\sin(u) - \mathrm{con}(\bar{f})u}{u} \right\|_C \right) \\
&= 0.
\end{aligned}$$

□

### 4.2. $\preccurlyeq_{\mathcal{P}_n^+}$-Optimization in Banach Space *C*[0, 1]

Let $n$ be a positive integer. Let $\mathcal{P}_n$ be the set of all polynomials with degrees less than or equal to $n$ and let $\mathcal{P}_n^+$ be the subset of $\mathcal{P}_n$ with nonnegative coefficients. Then $\mathcal{P}_n^+$ is a nontrivial closed convex and pointed cone in $C[0, 1]$. And therefore, $\mathcal{P}_n^+$ induces a partial order on $C[0, 1]$, which is denoted by $\preccurlyeq_{\mathcal{P}_n^+}$. Then, by definition, we have

$$f \preccurlyeq_{\mathcal{P}_n^+} g \quad \Leftrightarrow \quad g - f \in \mathcal{P}_n^+, \text{ for } f, g \in C[0,1]. \tag{4.6}$$

So, by definition (4.6), we have that $f \preccurlyeq_{\mathcal{P}_n^+} g$ if and only if $g - f$ is a polynomial with degree less than or equal to $n$ and with nonnegative coefficients. Let $\mathbb{Z}$ denote the set of all integers

**Proposition 4.4**. *Let* $S: C[0,1] \to C[0,1]$ *be the single-valued mapping defined in Theorem* 4.3. *Then*

$$\nabla S^0 = dS^0 = \{n\pi + \tfrac{\pi}{2} : n \in \mathbb{Z}\}$$

*and* $$E(S, \preccurlyeq_{\mathcal{P}_n^+}) = \emptyset. \tag{4.7}$$

*Here for each* $n \in \mathbb{Z}$, $n\pi + \frac{\pi}{2}$ *defines a constant function on* [0, 1].

*Proof.* Let $\theta_{C,C}$ be the zero continuous and linear mapping on $C[0,1]$. For given $\bar{f} \in C[0,1]$, by (4.5) in Theorem 4.3, we have that

$$\nabla S(\bar{f}) = \theta_{C,C} \quad \Leftrightarrow \quad \text{con}(\bar{f}(t)) = 0, \text{ for each } t \in [0, 1].$$

This implies that

$$\nabla S(\bar{f}) = \theta_{C,C} \quad \Leftrightarrow \quad \bar{f} = n\pi + \frac{\pi}{2}, \text{ for some } n \in \mathbb{Z}.$$

This proves the first part in this theorem. For any given $n\pi + \frac{\pi}{2} \in \nabla S^0$, we have $S\left(n\pi + \frac{\pi}{2}\right) = \pm 1$. In general, for $g \in C[0, 1]$ with $g$ is not a constant, then $S(g) = \sin(g)$ is not a polynomial with degree less than or equal to $n$, which implies that $\sin(g) \notin \pm\mathcal{P}_n^+$. Hence, if for $g$ is not a constant function defined on [0, 1], then

$$S\left(n\pi + \frac{\pi}{2}\right) \text{ and } S(g) \text{ are not } \preccurlyeq_{\mathcal{P}_n^+}\text{-comparable.} \tag{4.8}$$

By the first part of (4.6) and part (i) in Corollary 2.5, (4.8) implies that $E(S, \preccurlyeq_{\mathcal{P}_n^+}) = \emptyset$. □

### 4.3. $\preccurlyeq_{C^+}$-Optimization in Banach Space *C*[0, 1]

It is clear that the ordered extrema of a given mapping depend on the partial order. In this subsection, we consider the Banach space *C*[0, 1] equipped with a partial order, which is different from the partial order $\preccurlyeq_{\mathcal{P}_n^+}$ on *C*[0, 1] studied in the previous subsection. We will find that with the same given single-valued mapping $S: C[0,1] \to C[0,1]$ defined in Theorem 4.3, its ordered extrema will be different with respect to different partial orders equipped on $C[0,1]$.

Let $C^+$ be the subset of *C*[0, 1] defined by

$$C^+ = \{f \in C[0,1]: f(t) \geq 0, \text{ for any } t \in [0,1]\}$$

One sees that $C^+$ is a nontrivial closed convex and pointed cone in *C*[0, 1], which induces a partial order in *C*[0, 1] denoted by $\preccurlyeq_{C^+}$. More precisely, for $f, g \in C[0, 1]$,

$$f \preccurlyeq_{C^+} g \quad \Leftrightarrow \quad g - f \in C^+. \tag{4.9}$$

So, by definitions of (4.9) and $C^+$, we have that $f \preccurlyeq_{C^+} g$ if and only if $g - f$ is a nonnegative function defined on $[0, 1]$. It is clearly to see that $\mathcal{P}_n^+ \subsetneqq C^+$. This immediately implies that, for $f, g \in C[0, 1]$,

$$f \preccurlyeq_{\mathcal{P}_n^+} g \quad \Rightarrow \quad f \preccurlyeq_{C^+} g,$$

and $$f \preccurlyeq_{C^+} g \quad \nRightarrow \quad f \preccurlyeq_{\mathcal{P}_n^+} g \tag{4.10}$$

**Proposition 4.5**. *Let* $S: C[0,1] \to C[0,1]$ *be the single-valued mapping defined in Theorem* 4.3. *Then*

(i) $E(S, \preccurlyeq_{C^+}) = \nabla S^0 = dS^0 = \{n\pi + \frac{\pi}{2} : n \in \mathbb{Z}\};$ (4.11)

(ii) $E^+(S, \preccurlyeq_{C^+}) = \{2n\pi + \frac{\pi}{2} : n \in \mathbb{Z}\};$

(iii) $E^-(S, \preccurlyeq_{C^+}) = \{(2n+1)\pi + \frac{\pi}{2} : n \in \mathbb{Z}\}$.

*Proof.* Proof of (ii). For any $n \in \mathbb{Z}$ and for any $f \in C[0, 1]$, we have

$$S(2n\pi + \tfrac{\pi}{2}) \equiv 1 \geq \sin(f(t)) = S(f)(t), \text{ for any } t \in [0, 1].$$

By (4.9), this implies that

$$S(2n\pi + \tfrac{\pi}{2}) \succcurlyeq_{C^+} S(f), \text{ for any } f \in C[0, 1]. \tag{4.12}$$

Hence, when $2n\pi + \frac{\pi}{2}$ is considered as a constant function defined on [0, 1], by (4.12), we have

$$2n\pi + \tfrac{\pi}{2} \in E^+(S, \preccurlyeq_{C^+}), \text{ for any } n \in \mathbb{Z}. \tag{4.13}$$

Similarly to the proof of (4.13), we can prove

$$(2n+1)\pi + \tfrac{\pi}{2} \in E^-(S, \preccurlyeq_{C^+}), \text{ for any } n \in \mathbb{Z}. \tag{4.14}$$

By (4.13) and (4.14), we obtain that

$$n\pi + \tfrac{\pi}{2} \in E(S, \preccurlyeq_{C^+}), \text{ for any } n \in \mathbb{Z}. \tag{4.15}$$

By $\nabla S^0 = dS^0 = \{n\pi + \frac{\pi}{2} : n \in \mathbb{Z}\}$ in Proposition 4.4, by Corollary 2.5 and by (4.15), part (i) of this proposition is proved. By part (i) and (4.13) and (4.14), parts (ii) and (iii) are proved, respectively. □

## 5. Gâteaux and Fréchet Derivatives of Mappings from $C[0, 1]$ to $L_p[0, 1]$

Let $(L_p[0, 1], \|\cdot\|_{L_p})$, with $1 < p < \infty$, be the standard Banach space of real valued function defined on [0.1] with null element $\theta_{L_p}$. Let $C[0, 1]$ be the Banach spaces studied in section 4. In this section, we investigate the explicit formulas of Gâteaux and Fréchet derivatives of some mappings from $C[0, 1]$ to $L_p[0, 1]$. We will see that the corresponding results are same with the mappings from $C[0, 1]$ to $C[0, 1]$ studied in the previous section.

**Theorem 5.1**. *Let m be a positive integer and let $R^m: C[0,1] \to L_p[0,1]$ be the $m^{th}$ power operator from $C[0, 1]$ to $L_p[0, 1]$ defined, for each $f \in C[0, 1]$, by*

$$R^m(f) = f^m \in L_p[0, 1].$$

*Then, $R^m$ is Fréchet differentiable on $C[0,1]$ such that for any given $\bar{f} \in C[0,1]$, the Fréchet derivative $\nabla R^m(\bar{f}): C[0,1] \to L_p[0,1]$ satisfies that, for any $g \in C[0,1]$,*

$$\nabla R^m(\bar{f})(g) = m\bar{f}^{m-1}g \in L_p[0, 1]. \tag{5.1}$$

*Proof.* Let $\bar{f} \in C[0,1]$ be arbitrarily given. Then, by the proof of Theorem 4.1, we estimate

$$\lim_{u \xrightarrow{C[0,1]} \theta_C} \left\| \frac{R^m(\bar{f}+u) - R^m(\bar{f}) - m\bar{f}^{m-1}u}{\|u\|_C} \right\|_{L_p}$$

$$= \lim_{u \xrightarrow{C[0,1]} \theta_C} \frac{\left\|(\bar{f}+u)^m - \bar{f}^m - m\bar{f}^{m-1}u\right\|_{L_p}}{\|u\|_C}$$

$$\leq \lim_{u \xrightarrow{C[0,1]} \theta_C} \frac{\left\|(\bar{f}+u)^m - \bar{f}^m - m\bar{f}^{m-1}u\right\|_{C}}{\|u\|_C}$$

$$= 0. \qquad \square$$

Notice that the polynomial type operator $P_m$ defined by (4.2) is also an operator from $C[0, 1]$ to $L_p[0,1]$. Then, similarly to Corollary 4.2, we have the following corollary of Theorem 5.1.

**Corollary 5.2**. *Let m be a positive integer and let* $P_m$: $C[0, 1] \to L_p[0,1]$ *be the polynomial type operator defined by* (4.2). *Then,* $P_m$ *has the following properties.*

(i) $P_m$ *is Fréchet differentiable on* $C[0, 1]$ *such that for any given* $\bar{f} \in C[0,1]$, *the Fréchet derivative* $\nabla P_m(\bar{f})$: $C[0,1] \to L_p[0,1]$ *satisfies that, for any* $g \in C[0,1]$,

$$\nabla P_m(\bar{f})(g) = \sum_{i=1}^{m} i a_i \bar{f}^{i-1} g \in L_p[0,1]. \tag{5.2}$$

(ii) $P_m$ *is Gâteaux differentiable on* $C[0, 1]$ *such that for any given* $\bar{f} \in C[0, 1]$, *the Gâteaux derivative* $(P_m)'(\bar{f})$ *satisfies that, for any* $g \in C[0, 1]$,

$$(P_m)'(\bar{f})(g) = \nabla P_m(\bar{f})(g) = \sum_{i=1}^{m} i a_i \bar{f}^{i-1} g \in L_p[0,1]. \tag{5.3}$$

**Theorem 5.3**. *Similarly to Theorem* 4.3, *we define a single-valued mapping* $T: C[0,1] \to L_p[0,1]$ *by*

$$T(f) = \sin(f), \text{ for each } f \in C[0, 1].$$

*Then,* $T$ *is Fréchet differentiable on* $C[0,1]$ *such that for any given* $\bar{f} \in C[0,1]$, *the Fréchet derivative* $\nabla T(\bar{f})$: $C[0, 1] \to L_p[0,1]$ *satisfies that*

$$\nabla T(\bar{f})(g) = \text{con}(\bar{f})g, \text{ for any } g \in C[0, 1].$$

*And therefore,* $T$ *is Gâteaux differentiable on* $C[0, 1]$ *and the Gâteaux derivative* $T'(\bar{f})$ *satisfies that*

$$T'(\bar{f})(g) = \nabla T(\bar{f})(g) = \text{con}(\bar{f})g, \text{ for any } g \in C[0,1]. \tag{5.5}$$

*Proof.* Let $\bar{f} \in C[0, 1]$ be arbitrarily given. Then, we calculate

$$\lim_{u \xrightarrow{C[0,1]} \theta_C} \left\| \frac{\sin(\bar{f}+u) - \sin(\bar{f}) - \text{con}(\bar{f})u}{\|u\|_C} \right\|_{L_p}$$

$$= \lim_{u \xrightarrow{C[0,1]} \theta_C} \left\| \frac{\sin(\bar{f})\text{con(u)} + \text{con}(\bar{f})\sin\text{(u)} - \sin(\bar{f}) - \text{con}(\bar{f})u}{\|u\|_C} \right\|_{L_p}$$

$$\leq \lim_{u \xrightarrow{C[0,1]} \theta_C} \left( \left\| \frac{\sin(\bar{f})\text{con(u)} - \sin(\bar{f})}{\|u\|_C} \right\|_{L_p} + \left\| \frac{\text{con}(\bar{f})\sin\text{(u)} - \text{con}(\bar{f})u}{\|u\|_C} \right\|_{L_p} \right)$$

$$\leq \lim_{u \xrightarrow{C[0,1]} \theta_C} \left( \left\| \frac{\sin(\bar{f})\operatorname{con}(\mathrm{u}) - \sin(\bar{f})}{\|u\|_C} \right\|_C + \left\| \frac{\operatorname{con}(\bar{f})\sin(\mathrm{u}) - \operatorname{con}(\bar{f})u}{\|u\|_C} \right\|_C \right)$$

$$\leq \lim_{u \xrightarrow{C[0,1]} \theta_C} \left( \left\| \frac{\sin(\bar{f})\operatorname{con}(\mathrm{u}) - \sin(\bar{f})}{u} \right\|_C + \left\| \frac{\operatorname{con}(\bar{f})\sin(\mathrm{u}) - \operatorname{con}(\bar{f})u}{u} \right\|_C \right)$$

$$= 0.$$ □

Let $L_p^+$ be the positive cone in $(L_p, \|\cdot\|_p)$ that is defined by

$$L_p^+ = \{f \in L_p : f(t) \geq 0, \text{for ever } t \in [0,1]\}.$$

It is clear that $L_p^+$ is a nontrivial closed convex and pointed cone in $L_p$. Let $\preccurlyeq_{L_p^+}$ be the partial order on $L_p$ induced by $L_p^+$. Similarly to Proposition 4.5, we have the following results for the mapping $T$.

**Proposition 5.4**. *Let* $T: C[0,1] \to C[0,1]$ *be the single-valued mapping defined in Theorem* 5.3. *Then*

(i) $E(T, \preccurlyeq_{L_p^+}) = \nabla T^0 = dT^0 = \{n\pi + \frac{\pi}{2} : n \in \mathbb{Z}\}$;

(ii) $E^+(T, \preccurlyeq_{L_p^+}) = \{2n\pi + \frac{\pi}{2} : n \in \mathbb{Z}\}$;

(iii) $E^-(T, \preccurlyeq_{L_p^+}) = \{(2n+1)\pi + \frac{\pi}{2} : n \in \mathbb{Z}\}$.

*Proof.* The proof of this proposition is similar to the proof of Proposition 4.5. □

**Author's declaration.** The author declares that he has no conflict of interest, and the manuscript has no associated data. This work was not supported by external funding.